\documentclass[11pt]{article}
\input amssym.def
\usepackage{graphicx}
\usepackage[latin1]{inputenc}
\def\sqw{\hbox{\rlap{\leavevmode\raise.3ex\hbox{$\sqcap$}}$%
\sqcup$}}

   \def\bz{\hbox{\it Z\hskip -4pt Z}}     
\def\demo{\noindent{\bf Proof \ }} \newtheorem{theorem}{Theorem} 
\newtheorem{lemma}{Lemma} \newtheorem{proposition}{Proposition} 
 \newtheorem{example}{Example} 
\newtheorem{remark}{Remark}  
\newtheorem{corollary}{Corollary}

\def\dim{{\rm \ dim}\,}
\def\projdim{{\rm \ projdim}\,}
\def\proj{{\rm \ proj}\,}

\def\cd{{\rm \ cd} }
\def\char{{\rm \ char} }

\def\pd{{\rm \ projdim}\  }

\def\card{{\rm \ card\,} }
\def\rad{{\rm \ rad\,} }
\def\ara{{\rm \ ara\,} }

\newcommand{\la}{\langle}
\newcommand{\ra}{\rangle}
\def\bv{\hbox{\bf V}}

\def\bz{\hbox{\it Z\hskip -4pt Z}}

\def\bp{\hbox{\it I\hskip -2pt P}}

\newcommand{\Fcal}{{\mathcal F}}

\newcommand{\Lcal}{{\mathcal L}}
\newcommand{\Ical}{{\mathcal I}}
\newcommand{\Jcal}{{\mathcal J}}
\newcommand{\Bcal}{{\mathcal B}}
\newcommand{\Acal}{{\mathcal A}}
\newcommand{\Kcal}{{\mathcal K}}
\newcommand{\Ccal}{{\mathcal C}}

\newcommand{\Scal}{{\mathcal S}}
\newcommand{\Pcal}{{\mathcal P}}

\newcommand{\Mcal}{{\mathcal M}}
\newcommand{\Qcal}{{\mathcal Q}}

\newcommand{\Dcal}{{\mathcal D}}

\def\demo{\noindent{\bf Proof \ }} 

\begin{document}

\begin{center}
\uppercase{{\bf  Equations of  2-linear ideals and arithmetical rank }}
\end{center}
\advance\baselineskip-3pt
\vspace{2\baselineskip}
\begin{center}
{{\sc Marcel
Morales}\\
{\small Universit\'e de Grenoble I, Institut Fourier, 
UMR 5582, B.P.74,\\
38402 Saint-Martin D'H\`eres Cedex,\\
and IUFM de Lyon, 5 rue Anselme,\\ 69317 Lyon Cedex (FRANCE)}\\
 }
\vspace{\baselineskip}
\end{center}

{\small \sc Abstract.}{\footnote { first version  February 2008}}  In this paper we consider reduced homogeneous ideals $\Jcal\subset S$ of a polynomial ring $S$, having a 2-linear resolution. 
\begin{enumerate}
\item We study   systems of generators of $\Jcal\subset S$.
\item We compute the arithmetical rank for a large class of projective curves having a 2-linear resolution.
\item  We show that the fiber cone $\proj \Fcal(I_{\Lcal})$ of  a lattice ideal $I_{\Lcal}$ of codimension two is a set theoretical complete intersection. 
\end{enumerate}
 
\section{Introduction} In this paper we work on reduced algebraic sets on the projective space $X\subset \bp^n$. At the end of the XIX century,  it was known by the italians geometers that if $X$ is irreducible then the degree of $X$ is always greater than the codimension of $X$ +1. They classified the varieties with minimal degree by geometric means. 
\begin{theorem}(Bertini, Castelnuovo, Del Pezzo)  Let $X\subset \bp^n$ be any irreducible nondegenerate variety of dimension $k$ having (minimal) degree 
$n-k+1$. Then $X$ is either
\begin{enumerate}
\item A quadric hypersurface; 
\item a cone over the Veronese surface in $\bp^5$; or 
\item a rational normal scroll.
\end{enumerate}
\end{theorem}
Later Jo Harris has given a complete  proof of the above results, in particular we  quote the following result from \cite{ha}, page 108  :

\noindent The ideal defining a  rational normal scroll is generated by the $2\times 2$ minors of the following  scroll matrix :
$$\Bcal=\pmatrix{L_{1,0}&L_{1,1}&...&L_{1,c_1}&\ ...&\ L_{r,0}&L_{r,1}&...&L_{r,c_r}\cr L_{1,1}&L_{1,2}&...&L_{1,c_1+1}&\ ...&\ L_{r,1}&L_{r,2}&...&L_{r,c_r+1} \cr},$$ 
where $L_{1,0},L_{1,2},...,L_{1,c_1+1},...,L_{r,1},L_{r,2},...,L_{r,c_r+1}$ are all linearly independent forms.

\noindent We consider this matrix as formed by blocks 
$$\Bcal_i=\pmatrix{L_{i,0}&L_{i,1}&...&L_{i,c_i}\cr L_{i,1}&L_{i,2}&...&L_{i,c_i+1}\cr    }.$$
If $c_i=0$ we will say that $\Bcal_i$ is a generic block, and if $c_i>0$  that $\Bcal_i$ is a non generic block.

It then follows that in the irreducible case, Harris has described a system of generators for the ideal of a variety having minimal degree. We can also consider the problem to describe set theoretically the varieties of minimal degree, and in particular determine the arithmetical rank, that is the minimal number of elements wich generates up to the radical the ideal of a variety having minimal degree. There are very few results in this direction we cite the followings:
 
\begin{theorem}(Verdi) \cite{ve}Let consider a rational normal scroll $\Ccal$  which defines a projective curve, its defining ideal $\Mcal$ is generated by the  the $2\times 2$ minors of the scroll matrix with only one non generic block: 
$$\Bcal=\pmatrix{L_{0}&L_{1}&...&L_{c}\cr L_{1}&L_{2}&...&L_{c+1}\cr    }.$$
Where $c>0$. For $j=1,...,c$ let  
$$F_j:= \sum_{k=0}^{j} (-1)^k {{j}\choose k} (L_{j+1} )^{j-k}L_k (L_j)^k$$
then $$\Mcal=\rad(F_1,...,F_{c}).$$ Since the codimension of $\Ccal$  is $c$ this implies that $\Ccal$ (or $\Mcal$) is a set theoretical complete intersection. 
\end{theorem}
On the other hand we have the following result which follows from the works of Hochster, Bruns-Vetta-Schwanzl:
\begin{theorem} \label{bruns} Let $\Bcal$ be a scroll matrix with $r\geq 2$ columns and only generic blocks, then 
\begin{itemize}
\item If $\char K=0,$ $cd (\Mcal)=\ara (\Mcal)=2r-3; \pd (S/\Mcal)=r-1$.
\item If $\char K=p>0$, then  $\ara (\Mcal)=2r-3$, and 
since $R/\Mcal$ is a Cohen-Macaulay ring, we have  that $\pd (S/\Mcal)=\cd (\Mcal)=r-1$.
 
\end{itemize}
 
\end{theorem}
 When the algebraic set is not irreducible, Xambo \cite{x} has described from the geometric point of view the algebraic sets beeing connected in codimension one and having minimal degree. The problem of describing generators of their ideals is much more complicated (see for example \cite{cep}), and a first work in this direction was done by Barile and Morales \cite{bm2}, \cite{bm4}. Later in \cite{eghp} it was introduced and developped the notion of linearly joined varieties. In my recent work \cite{mo-stci}, I study the case of linearly joined linear spaces, there I have computed the projective dimension, the cohomological dimension and the arithmetical rank. 
 
The purpose of this paper (see section 2), is to describe more precisely the generators  of ideals  having a 2-linear resolution. For ideals having 2-linear resolution we refer the  reader to \cite{eg}, \cite{bm2}, \cite{bm4}, \cite{eghp}, \cite{mo-stci}.

In section 3 we consider a large class of ideals  having a 2-linear resolution (in fact projective curves) and we compute its arithmetical rank. This work extends all the works in this direction in the litterature, in particular those very interesting of Margherita Barile \cite{b1}.

Our results apply to the special case where $X$ is the fiber cone $\proj \Fcal(I_{\Lcal})$ of  a lattice ideal $I_{\Lcal}$ of codimension two. Let recall that by the succesively works \cite{g}, \cite{gms1}, \cite{gms2},\cite{bm1} and \cite{hm}, the fiber cone is a reduced arithmetically Cohen-Macaulay algebraic set of minimal degree. See section 4.
\section{Equations of  2-linear ideals.}
We recall the following result which follows from \cite{bm4} and \cite{eghp}
\begin{theorem}{\label{2-linear}} The following conditions are equivalent:
\begin{enumerate}
\item  The reduced ideal   $\Jcal\subset S:=K[\bv]$ has 
 {\it 2-linear resolution}, 
\item  Let $\Jcal=\Jcal_1\cap ...\cap \Jcal_l$ be the primary decomposition of $\Jcal\subset S:=K[\bv]$. For all $i=1,...,l$, there exist   sublinear spaces  $\Dcal_i,\Pcal_i \subset \bv$, with $\Dcal_l=0,\Pcal_1=0$,
 and ideals  $\Mcal_i \subset K[\bv]$ such that 
\begin{itemize}
\item a) for all $i=1,...,l$, $\Jcal_i=(\Mcal_i, \Qcal_i)$, $\Mcal_i  $  is the ideal of the $2\times 2$ minors of a scroll matrix,
\item b) $\Qcal_i =
\Dcal_{i } \oplus  \Pcal_{i },$   
\item c) $\Dcal_i =
\Delta_{i+1 } \oplus ...\oplus \Delta_{l},$ 
\item d) $\Mcal_i\subseteq (\Delta _{i })$ for all $i=2,...,l$,
\item e) $\Mcal_i\subseteq (\Pcal_{j })$ for all $1\leq i<j\leq l$,
\item f) $\bigcap_{j=1 }^{ k-1}(\Qcal_j)\subseteq (\Pcal_k,\Dcal_{k-1 })$ for all $k=2,...,l$.
\end{itemize}
\end{enumerate}
\end{theorem} 

In the  following result we will try to have more information on the ideals $\Mcal_i$ :
\begin{lemma} Let $\Bcal_i$ be a non generic block of a scroll matrix, with  $c_i+1$ columns,
$$\Bcal_i=\pmatrix{L_{i,0}&L_{i,1}&...&L_{i,c_i}\cr L_{i,1}&L_{i,2}&...&L_{i,c_i+1}\cr    },$$
 suppose that the ideal $\Mcal_i$ generated by all the $2\times 2$ minors of $\Bcal_i$, is contained in a linear ideal $(\Delta )$. Then for all $i=1,...,r$ we have either
\begin{enumerate}
\item $\{L_{i,0},L_{i,1},...,L_{i,c_i} \}\subset \Delta $ or $\{L_{i,1},L_{i,2},...,L_{i,c_i+1} \}\subset \Delta $
\item there exists a linear form $H\notin \Delta $ and a non zero  constant  $\alpha_i$ such that $\{L_{i,0}-H,L_{i,2}-\alpha _iH,...,L_{i,c_i+1}-\alpha _i^{c_i+1}H \}\subset \Delta $
\end{enumerate}

\end{lemma}
\demo Let $U\subset \bv$ be a subvector space  such that  
$\bv=U \oplus\Delta $, we can write 
$L_{i,j}=L'_{i,j}+ H_{i,j}$ with 
$L'_{i,j}\in \Delta  , H_{i,j}\in U.$ Since for all $j\geq 0$,
 $L_{i,j} L_{i,j+2}-L_{i,j+1}^2\in (\Delta )$, we have 
$H_{i,j} H_{i,j+2}-H_{i,j+1}^2=0$ 
We have to consider two cases:
\begin{enumerate}
\item $H_{i,j}=0$ for some $1\leq j\leq  c_i$ this implies $H_{i,1}=...= H_{i,c_i}=0$, and from the relation $L_{i,0} L_{i,c_i+1}-L_{i,1}L_{i,c_i}\in (\Delta )$ we have that $H_{i,0} H_{i,c_i+1}=0$ so that we have either $H_{i,0} =0$ or $H_{i,c_i+1}=0$. 
\item $H_{i,j}\not=0$ for all $1\leq j\leq  c_i$, the relations  $H_{i,j} H_{i,j+2}-H_{i,j+1}^2=0$ imply 
that there exist a non zero  constant  $\alpha_i$ such that $H_{i,j}=\alpha _i^jH_{i,0},$ and so 
$L_{i,j}=L'_{i,j}+ \alpha _i^jH_{i,0}$  with 
$L'_{i,j}\in \Delta  , H_{i,0}\notin \Delta .$
\end{enumerate}
\begin{lemma} \label{generic} Let $$\Bcal=\pmatrix{L_{1,0} &L_{2,0} &...&L_{r,0}\cr L_{1,1} &L_{2,1} &...&L_{r,1}\cr}$$ be a $2\times r$ generic matrix (with $r\geq 2$),  suppose that the ideal $\Mcal$ generated by all the $2\times 2$ minors of $\Bcal$, is contained in a linear ideal $(\Delta )$. Then  we have either
\begin{enumerate}
\item  Up to permutation of the lines of $\Bcal$, there exists some $i$ such that $L_{i,0}\in \Delta , L_{i,1}\notin \Delta $; in this case  the elements of the first line of $\Bcal$ are in $ \Delta $,   
\item  there exists some $i$ such that $L_{i,0},L_{i,1}\in \Delta $; in this case the matrix obtained from $\Bcal$ by deleting the column $i$, $\Bcal'$ is still a generic matrix with $r-1$ columns, and if $r-1\geq 2$ then we apply recursively this lemma to $\Bcal'$,
\item for any $i=1,...,r$ there exists a nonzero linear form $H_{i}\notin \Delta $, and a non zero  constant  $\alpha$ such that $\{L_{{i},0}-H_{i},L_{{i},1}-\alpha H_{i} \}\subset \Delta $, 
\item for any $i=1,...,r$ there exists  linear forms $H_{1},H_2 \notin \Delta$, and a non zero  constant  $\alpha_i$ such that $\{L_{{i},0}-\alpha_i H_{1},L_{{i},1}- \alpha_i H_{2} \}\subset \Delta .$ 
\end{enumerate}

\end{lemma}
\demo 
Let $U$ be a subvector space of $\bv$ such that $\bv=U\oplus \Delta $, so that we can write $L_{{i},j}=L'_{{i},j} + H_{{i},j}$, with $ L'_{{i},j} \in \Delta,  H_{{i},j}\in U$. Since $L_{{1},0}L_{{i},1}-L_{{1},1}L_{{i},0}\in (\Delta ) $ we have that 
$H_{{1},0}H_{{i},1}-H_{{1},1}H_{{i},0}=0 $, we have two choices:
\begin{enumerate}
\item $H_{{i},j}=0$ for some $i,j$, then up to a permutation of the lines and columns  we can assume that $H_{{1},0}=0$. If $H_{{1},1}\not=0$ the above relation implies that $H_{{i},0}=0$ and then all the entries of the  first line of $\Bcal$ are in $\Delta $.
 If both $H_{{1},0}=0, H_{{1},1}=0$, then  $L_{1,0},L_{1,1}\in \Delta , $ in this case the matrix obtained from $\Bcal$ by deleting the first column, $\Bcal'$ is still a generic matrix with $r-1$ columns,
\item $H_{{i},j}\not=0$ for all $i,j$, then 
\begin{itemize}
\item there exist  a nonzero constant  $\alpha $ such that  $H_{{1},1}= \alpha H_{{0},0}$, and this implies  $H_{{i},1}= \alpha H_{{i},0}$, for all $i$
\item there exist  nonzero  constants  $\alpha_i $ such that $H_{{i},0}= \alpha_i H_{{0},0}$, and this implies  $H_{{i},1}= \alpha_i H_{{0},1}$, for all $i$.
\end{itemize}

\end{enumerate}
\begin{proposition} Let $\Bcal$ be a scroll matrix, suppose that the ideal $\Mcal$ generated by  the $2\times 2$ minors of $\Bcal$  is contained in a linear ideal $(\Delta ) $. We have either:
\begin{enumerate}
\item  All the blocks of $\Bcal$ are generic,  so Lemma \ref{generic} applies.
\item The entries of exactly one  line of $\Bcal$ are in $\Delta $.
\item there exist one  block $\Bcal_i$ of $\Bcal$ which is non generic, and all the entries of $\Bcal_i$ are in $\Delta $  then the matrix $\Bcal'$ obtained from  $\Bcal$ by deleting $\Bcal_i$ is a scroll matrix with less columns than $\Bcal$ and we must apply recursively this proposition to $\Bcal'$.
\item there exist one  $\Bcal$ which is non generic, and any entry of a non generic block $\Bcal_i$ is not in $\Delta$.  In this case  there exists a non zero  constant  $\alpha$, and 
\begin{itemize}
\item for any non generic block $\Bcal_i$ there exists a linear form  $H_{i}\notin \Delta $,  such that $\{L_{{i},0}-H_{i},L_{{i},1}-\alpha H_{i},...,L_{{i},c_{i}+1}-\alpha ^{c_{i}+1}H_{i} \}\subset \Delta $.
\item for any generic block $\Bcal_i$, there exists a linear form $H_{i} $, null or $H_{i}\notin \Delta $,  such that $\{L_{{i},0}-H_{i},L_{{i},1}-\alpha H_{i} \}\subset \Delta $.
\end{itemize}
\end{enumerate}
\end{proposition}
\demo The proof follows easily from the above lemmas.
\section{Arithmetical rank of some 2-linear ideals.}
From now on, we suppose that $\Jcal$ is a reduced ideal having a $2-$linear resolution, we know  the primary decomposition 
$$ \Jcal=(\Mcal_1, \Qcal_1)\cap ... \cap (\Mcal_l, \Qcal_l)$$
where $\Qcal_i =
\Dcal_{i } \oplus  \Pcal_{i },$   
 $\Dcal_i =
\Delta_{2 } \oplus ...\oplus \Delta_{i +1},$ 
 $\Mcal_i\subseteq (\Delta _{i })$ for all $i=2,...,l$,
 $\Mcal_i\subseteq (\Pcal_{j })$ for all $i=1,...,l-1$ and $j=i+1,...,l$,
 We recall the following theorem from \cite{mo-stci}
 \begin{theorem}Let $\Jcal$ be a reduced ideal having a $2-$linear resolution, with the above notations 
 \begin{enumerate}
\item Let $\Qcal:=\Qcal_1\cap ... \cap  \Qcal_l$, then  $\pd (S/\Jcal)=\pd (S/\Qcal)=\max _{i=2,...,l}\{\dim_K (\Pcal_i+\Dcal_{i-1})-1\}.$
\item With the asumption that for all $i$, $\Mcal_i$ is a set theoretical complete intersection, we have 
$$\cd \Jcal= \pd (S/\Jcal)=\max _{i=2,...,l}\{\dim_K (\Pcal_i+\Dcal_{i-1})-1\}.$$

\end{enumerate}

\end{theorem}

\begin{theorem}\label{ara}
Let $\Jcal$ be a reduced ideal having a $2-$linear resolution, assume that for all $i$  the matrix $\Bcal_i$ is either the matrix null or has at most one  block and one line of $\Mcal_i$ is included in $ (\Delta _{i })$ for all $i=2,...,l$, and one line of $\Mcal_i$ is included in $ (\Pcal_{j })$ for all $i=1,...,l-1$ and $j=i+1,...,l$. Then 
$$ara(\Jcal)=cd(\Jcal)=\projdim (S/\Jcal).$$
\end{theorem} 
\demo Let 
$$\Bcal_i=\pmatrix{L_{i,0}&L_{i,1}&...&L_{i,c_i}\cr L_{i,1}&L_{i,2}&...&L_{i,c_i+1}\cr    },$$ and $\Mcal_i$ the ideal generated by all the $2\times 2$ minors of $\Bcal_i$. We call $L_{i,0}, L_{i,c_i+1}$, "corner" entries of $\Bcal_i$ and the other "inner" entries of $\Bcal_i$.
 We know that $\Jcal$ is generated by all the ideals $\Mcal_1,...,\Mcal_l$ and   $\bigcup _{j=1 }^{ k} 
( \Delta_{j } \times \Pcal_j),$
where $(\Delta_{j } \times \Pcal_j)$ is the ideal generated by all the products $fg$, with 
$f\in\Delta_{j } , g\in \Pcal_j.$ We  have also the property that $ \Delta_{j } \times \Pcal_j\subset \Pcal_i$ when $j<i$.

First we study the intersection of the two primary components of $\Jcal$. We have that $(\Mcal_1, \Qcal_1)\cap (\Mcal_2, \Qcal_2)= (\Dcal_2,(\Mcal_1, \Mcal_2,\la \Delta_{2 }\ra \times \Pcal_2)$ Since $\Mcal_1 \subset \Pcal_2$, by our hypothesis we can suppose that the first line of $\Bcal_1$, $\Lcal_1:=\{L_{1,0},L_{1,1},...,L_{1,c_{1}} \}$ is included in $ \Pcal_2 $, and we can find a vector space  $\widetilde{\Pcal_2 }$ containing a corner element of $\Bcal_1$  such that $\Pcal_2 =\widetilde{\Pcal_2 } \oplus \widetilde{\Lcal_1 }$, where $\widetilde{\Lcal_1 }$ is the vector space generated by all the inner entries of $\Bcal_1$. By the same reasons 
we can find a vector space  $\widetilde{\Delta _2 }$ containing a corner element of $\Bcal_2$  such that $\Delta _2 =\widetilde{\Delta _2 } \oplus \widetilde{\Lcal_2 }$, where $\widetilde{\Lcal_2 }$ is the vector space generated by all the inner entries of $\Bcal_2$.

Suppose that $L_{1,0}\in \widetilde{\Pcal_2 }$, working modulo $\Mcal_1  $ we have the relations $L_{1,0} L_{1,2}-L_{1,1}^2=0$ so for any element $g\in \Delta _2$ we will have  
$(L_{1,0}g) (L_{1,2}g)=(L_{1,1}g)^2,$ so that $(L_{1,1}g)^2\in (\Mcal_1,(L_{1,0}g))$, since we have also 
$(L_{1,1}g) (L_{1,3}g)=(L_{1,2}g)^2,$ then $(L_{1,2}g)^2\in (\Mcal_1,(L_{1,1}g))$, by repeating this argument we have that there exists some powers $n_i$ such that $(L_{1,i}g)^{n_i}, \in (\Mcal_1,(L_{1,0}g))$. Let $L_{2,0}\in \widetilde{\Delta _2 }$ a "corner" entry of $\Bcal_2$, we use the relations in the matrix $\Bcal_2$,  then  for all "inner" element $L$ of  $\Bcal_2$ and every element $h\in \Pcal_2$ there exists some powers $m_i$ such that $(Lg)^{m_i}, \in (\Mcal_2,(L_{2,0}h))$. 

In other words, let $\widetilde{\Lcal_i}$ the vector space generated by the inner entries of the matrix $\Bcal_i$, We have found    vector spaces  $\widetilde{\Pcal_2 }, \widetilde{\Delta _2 }$, such that 
$$ \Pcal_2= \widetilde{\Pcal_2 }\oplus \widetilde{\Lcal_1} ,  \Delta _2= \widetilde{\Delta _2 }\oplus \widetilde{\Lcal_1}     ,$$ 
and if  $\Ical_2$ is the ideal  generated by  $\Mcal_1,\Mcal_2$ and   $\widetilde{\Delta_2 }\times \widetilde{\Pcal_2 } $ then up to radical the following ideals $(\Ical_2,\Dcal_2)$,$(\Mcal_1, \Qcal_1)\cap (\Mcal_2, \Qcal_2)$ are equal.
Let remark that $\dim \widetilde{\Pcal_2 }= \dim {\Pcal_2 }-c_1$,  $\dim \widetilde{\Delta _2 }= \dim {\Delta _2 }-c_2$.

We will prove the following statement:

For any $i\geq 2$, there exists a decomposition 
$$ \Pcal_i= \widetilde{\Pcal_i }\oplus \widetilde{\Lcal_1}\oplus...\oplus  \widetilde{\Lcal_{i-1}} ,  \Delta _i= \widetilde{\Delta _i }\oplus \widetilde{\Lcal_i}     ,$$ 
such that 
\begin{description}
\item [($H_{i,1}$)] for any $j<k\leq i$ we have $\widetilde{\Pcal_j }\times \widetilde{\Delta _j }\subset (\widetilde{\Pcal_k })$
\item [($H_{i,2}$)] Let $\Ical_i$ be the ideal  generated by  $\Mcal_1,...,\Mcal_i$ and   $\bigcup_{j=2}^{i}(\widetilde{\Delta_j }\times \widetilde{\Pcal_j }).$ Up to radical the ideals  $(\Ical_i,\Dcal_i)$, and $\bigcap_{j=2}^{i}(\Mcal_i, \Qcal_i)$ are equal. 
\end{description}
We suppose that this is true for some $i\geq 2$ and prove it for $i+1$.
Since $\Delta _j \times\Pcal_j \subset (\Pcal_{i+1})$ for any $j<i+1$, and since $(\Pcal_{i+1})$ is a prime ideal, for any $j<i+1$ we have either $\Delta _j  \subset \Pcal_{i+1}$ or $\Pcal_j \subset \Pcal_{i+1}$ (equality like vector spaces), we have two cases :
\begin{enumerate}
\item $\Delta _j  \subset \Pcal_{i+1}$ for all $j=2,...,i$
\item there exists $2\leq i_0\leq i$ such that $\Pcal_{i_0} \subset \Pcal_{i+1}$ and $\Delta _j  \subset \Pcal_{i+1}$ for all $j=i_0+1,...,i$
\end{enumerate}
Since $\Mcal_{i+1} \subset (\Delta _{i+1})$, then by our hypothesis a line of $\Bcal_{i+1} $ is included in $\Delta _{i+1}$, let 
$\widetilde{\Delta _{i+1} }$ be a subvector space containing the corner element in this line such that 
$ \Delta _{i+1}=\widetilde{\Delta _{i+1} }\oplus \Lcal_{i+1}.$ 

In the first case, it follows that  $\Pcal_{i+1}$ contains $\bigoplus_{j=2}^{i}\Delta _j$ and since $\Mcal_1 \subset (\Pcal_{i+1})$, by our hypothesis $\Pcal_{i+1}$ contains a line in the matrix $\Bcal_1$, formed by a corner element $L_1$ and the inner elements, $\widetilde{\Lcal_1}  \subset \Pcal_{i+1}$, so $\Pcal_{i+1}\supset (\bigoplus_{j=2}^{i}\widetilde{\Delta _j }\oplus KL_1) \oplus \bigoplus_{j=2}^{i}\widetilde{\Lcal_j }$. Let $\widetilde{\Pcal_{i+1}} $
a vector space containing $\bigoplus_{j=2}^{i}\widetilde{\Delta _j }\oplus KL_1$ such that 
$\Pcal_{i+1}= \widetilde{\Pcal_{i+1}}\oplus \bigoplus_{j=2}^{i}\widetilde{\Lcal_j }$. It is clear that ($H_{i+1,1}$), ($H_{i+1,2}$) are satisfied in this case.

In the second case, it follows that  $\Pcal_{i+1}$ contains $\bigoplus_{j=i_0+1}^{i}\Delta _j\oplus \Pcal_{i_0}$ and since $\Mcal_{i_0} \subset (\Pcal_{i+1})$, by our hypothesis $\Pcal_{i+1}$ contains a line in the matrix $\Bcal_{i_0}$, formed by a corner element $L_{i_0}$ and the inner elements, $\widetilde{\Lcal_{i_0}}  \subset \Pcal_{i+1}$, so $\Pcal_{i+1}\supset (\bigoplus_{j=i_0+1}^{i}\widetilde{\Delta _j }\oplus KL_{i_0})\oplus \widetilde{\Pcal_{i_0}}  \oplus \bigoplus_{j=2}^{i}\widetilde{\Lcal_j }$. Let $\widetilde{\Pcal_{i+1}} $
a vector space containing $\bigoplus_{j=2}^{i}\widetilde{\Delta _j }\oplus KL_1$ such that 
$\Pcal_{i+1}= \widetilde{\Pcal_{i+1}}\oplus \bigoplus_{j=2}^{i}\widetilde{\Lcal_j }$. It is clear that ($H_{i+1,1}$), ($H_{i+1,2}$) are satisfied in this case.

We have used the property that for any $k=1,...,l$ the entries of the matrix $\Bcal_k$ are linearly independent and also linearly independent with $\bigoplus_{j=k+1}^{l}\Delta _j\oplus \Pcal_{k}$.

As a conclusion, $\Jcal$ us up to the radical equal to the ideal $\Ical$ generated by $\Mcal_1,...,\Mcal_i$ and   $\bigcup_{j=2}^{l}(\widetilde{\Delta_j }\times \widetilde{\Pcal_j }).$ 
On the other hand the ideal $\Kcal:=\bigcup_{j=2}^{l}(\widetilde{\Delta_j }\times \widetilde{\Pcal_j })$ has a 2-resolution, with primary decomposition : $\Kcal= \bigcap_{j=1}^{l}(\widetilde{\Dcal_j }, \widetilde{\Pcal_j }  )$, where $\widetilde{\Dcal_j }:=\bigoplus_{k=j+1}^{l}\widetilde{\Delta_k }$ and so 
$\ara(\Kcal)=\max_{2\leq j\leq l} \{\dim \widetilde{\Dcal_{j-1} }+\dim \widetilde{\Pcal_j }\}-1$
by a quick computation we have that 
$$ \dim \widetilde{\Dcal_{j-1} }+\dim \widetilde{\Pcal_j }= \dim {\Dcal_{j-1} }+\dim {\Pcal_j }-(c_1+...+c_l)$$
so that 
$$\ara(\Jcal)\leq  \ara(\Kcal)+(c_1+...+c_l)= \max_{2\leq j\leq l} \{\dim {\Dcal_{j-1} }+\dim {\Pcal_j }\}-1=\projdim (S/\Jcal).$$
On the other hand by \cite{mo-stci} we have that $\cd \Jcal=\projdim (S/\Jcal)\leq \ara(\Jcal),$ the equality follows.
\begin{remark}
We can apply the method of the proof of the above theorem in order to get an upper bound for the arithmetical rank of any $2-$linear ideal in terms of the arithmetical rank of rational normal scrolls. Note that the problem to compute the arithmetical rank of rational normal scrolls is still open in general.
\end{remark} 
\begin{theorem}\label{ara1}
Let $\Jcal$ be a reduced ideal having a $2-$linear resolution, $$ \Jcal=(\Mcal_1, \Qcal_1)\cap ... \cap (\Mcal_l, \Qcal_l),$$ 
where $\Mcal_i$ is either the ideal null or the ideal generated by all the $2\times 2 $ minors of a scroll matrix $\Bcal^{(i)}$.  Assume that $\Bcal^{(i)}$ has $j_i$  non generic blocks, and the $k$-non-generic block  in $\Bcal^{(i)}$ has  $c_{i,k}+1$  columns.   Moreover assume that one line of $\Bcal^{(i)}$ is included in $ (\Delta _{i })$ for all $i=2,...,l$, and one line of $\Bcal^{(i)}$ is included in $ (\Pcal_{j })$ for all $1\leq i<j\leq l$.  Setting $\overline {c}_i=c_{i,1}+...+c_{i,j_i}$ we have that  
$$\ara(\Jcal)\leq \projdim (S/\Jcal)+ \sum_{i=1}^l(\ara \Mcal_i -\overline {c}_i).$$
\end{theorem} 
 \begin{example} 
Consider  $S=K[a,b,c,x,y,z,u,v,w]$ the   ring of polynomials, and 
$$\Jcal_1=(uv-w^2,a,b,c); \Jcal_2=(y,z,v,w,a,b); \Jcal_3=(x,z-u,v,w,c,b); \Jcal_4=(x-u,y-u,a,c,v,w).$$
The sequence $\Jcal_1, ...,\Jcal_4 $ satisfies the condition 2 of Theorem \ref{2-linear}, so $\Jcal:=\bigcap_{i=1}^{4}\Jcal_i$ has a 2-linear resolution, 
 $\projdim(S/\bigcap_{i=1}^{4}\Jcal_i)=6$ and 
$\Jcal= \rad( uv-w2,cb, ca+ab, cy+ax+b(x-u), cz+a(z-u)+b(y-u),cv+av+bv).$

\end{example}
 \begin{example} 
Consider  $S=K[a,b,c,x,y,z,u,v,w]$ the   ring of polynomials, and 
$$\Jcal_1=(a,b,c,x); \Jcal_2=(x(x-u)-c^2,a,b,y,z); \Jcal_3=(b,x,z-u,c); \Jcal_4=(x,y-u,a,c).$$
The sequence $\Jcal_1, ...,\Jcal_4 $ satisfies the condition 2 of Theorem \ref{2-linear}, so  $\Jcal:=\bigcap_{i=1}^{4}\Jcal_i$ has a 2-linear resolution, 
$ \Jcal$ is generated by the polynomial $x(x-u)-c^2$ and the elements in the following tableau:
$$ bc$$
$$bx\hskip 2cm ac$$
$$ba\hskip 2cm  ax \hskip 2cm cy$$
$$b(y-u)\hskip 0.6cm  a(z-u)\hskip 1cm  xy\hskip 2cm  cz$$
$$ \hskip 3.5cm  xz,$$
 $\projdim(S/\bigcap_{i=1}^{4}\Jcal_i)=5,$ up to the radical $ \Jcal$ is generated by the polynomial $x(x-u)-c^2$ and the elements in the following tableau: 
$$bx$$
$$ba\hskip 2cm  ax $$
$$b(y-u)\hskip 0.6cm  a(z-u)\hskip 1cm  xy$$
$$ \hskip 6cm  xz$$
thus $\ara \Jcal=5,$ and $\Jcal= \rad(x(x-u)-c^2, bx,ab+ax,b(y-u)+a(z-u)+xy,xz ).$

\end{example}
\begin{example} In this example we can compute the arithmetical rank by using a general theorem  due to Eisenbud and Evans.
Consider  $S=K[x_1,...,x_{2r},\Delta ]$ a  ring of polynomials where $\Delta$ is a set of variables, $K$ is a field of characteristic zero, let $0\leq \alpha\leq r, $
$$\Jcal_1=(\mid\matrix{x_1&...&x_{r}\cr x_{r+1}&...&x_{2r}\cr}\mid,\Delta ); \Jcal_2=(x_1,...,x_{2r-\alpha }); $$
and  $\Jcal:=\Jcal_1\cap \Jcal_2$, then $\Jcal$ has a 2-linear resolution.
We have the following exact sequence:
$$0\rightarrow S/\Jcal_1 \cap \Jcal_2\rightarrow S/\Jcal_1\oplus  S/\Jcal_2\rightarrow S/(\Jcal_1 + \Jcal_2)\rightarrow 0,$$
 which gives rise to the 
long exact sequence:
$$\rightarrow H_{\Jcal_1 \cap \Jcal_2}^{h-1}(S)\rightarrow H_{\Jcal_1 + \Jcal_2}^{h}(S)\rightarrow H_{\Jcal_1 }^{h}(S)
\oplus H_{\Jcal_2}^{h}(S)\rightarrow H_{\Jcal_1 \cap \Jcal_2}^{h}(S)\rightarrow H_{\Jcal_1 + \Jcal_2}^{h+1}(S)\rightarrow .$$ 
By using the Bruns-Vetta Schwanzl's theorem \ref{bruns}   and Proposition 2 of \cite{ha-mo} we have that $\cd \Jcal_1=(2r-3)+\card \Delta  $, on the other hand  $\cd \Jcal_2= 2r-\alpha $,$\cd (\Jcal_1 + \Jcal_2)=\cd (\Delta ,x_1,...,x_{2r-\alpha})=2r-\alpha+\card \Delta $, so if $0\leq \alpha\leq 2, $ we have   $\cd \Jcal=\dim S -\alpha -1$, by \cite{mo-stci} we get   
 $\projdim(S/\Jcal)=\dim S -\alpha-1=\cd \Jcal,$  
and by  \cite{ee}  $\ara \Jcal\leq \dim S -1$,  finally $\dim S -\alpha-1\leq \ara \Jcal\leq  \dim S -1, $ if $\alpha =0$ we have  $\ara \Jcal= \dim S -1.$
\end{example}
\begin{example}More generally, consider a ring of polynomials $S=K[\Delta,P ]$, where $\Delta ,P$ are disjoints sets of variables, $\Mcal\subset (P)$ be any ideal, let
$$\Jcal_1=(\Mcal,\Delta ); \Jcal_2=(P); \Jcal:=\Jcal_1\cap \Jcal_2,$$
then by the  argument developped in the above example we have that $\cd \Jcal=\dim S -1$,  and  by  \cite{ee}  $\ara \Jcal\leq \dim S -1$, so  $\ara \Jcal= \dim S -1$.   
\end{example} 
 We give now some simple open cases.
\begin{example} 
Consider  $S=K[a,b,c,d,e,f]$ the   ring of polynomials,  $\Bcal=\pmatrix{c-f&d-f\cr d-f&d+c-f\cr}$ and $F$ the determinant of $\Bcal$, then the ideal $\Jcal=\bigcap_{i=1}^{3}\Jcal_i$, where
$$\Jcal_1=(F,a,b); \Jcal_2=(b,c,d); \Jcal_3=(c,d,e); $$
has a 2-linear resolution, 
 $\cd (\Jcal)=\projdim(S/\Jcal)=3$ and 
$\Jcal= \rad( F,bc, ac+bd, ad+be),$ so $3\leq \ara (\Jcal)\leq 4$. We guess that $ \ara (\Jcal)= 4.$

\end{example}
\begin{example} 
Consider  $S=K[a,b,c,d,\Delta ]$ a   ring of polynomials,  $\Bcal=\pmatrix{a&c\cr b&d\cr}$, then the ideal $\Jcal=\Jcal_1\cap \Jcal_2$, where
$$\Jcal_1=(ad-bc,\Delta ); \Jcal_2=(b,d)   $$
has a 2-linear resolution, 
 $\cd (\Jcal)=\projdim(S/\Jcal)=\card(\Delta )+1$. In \cite{b1} it is proved that 
$ \ara (\Jcal)= \card(\Delta )+1.$

\end{example}
\begin{example} 
Consider  $S=K[a,b,c,d,e,f,g]$ the   ring of polynomials,  $\Bcal=\pmatrix{a&c\cr b&d\cr}$ and $F=ad-bc$, then the ideal $\Jcal=\bigcap_{i=1}^{3}\Jcal_i$, where
$$\Jcal_1=(F,e,f); \Jcal_2=(b,d,f); \Jcal_3=(b,d,g); $$
has a 2-linear resolution, 
 $\cd (\Jcal)=\projdim(S/\Jcal)=3$. We will prove that $\ara (\Jcal)=3.$
 Set $q_1:=aF+be; $
 and  $q'_1:=a^2q_1+bf, q'_2:=cF+de+fg; q'_3:=(ac-e)q_1+df$, we assert that
$\Jcal= \rad(q'_1,q'_2,q'_3 ),$ so $\ara (\Jcal)=3$. 
First we have that $dq'_1-bq'_3=q_1^2$ so $q_1\in \rad(q'_1,q'_2,q'_3 ),$ and it follows that $bf, df\in \rad(q'_1,q'_2,q'_3 ).$ We have also that $dq_1-bq'_2=F^2 -bfg$ which implies that $F\in \rad(q'_1,q'_2,q'_3 ),$ and then  $be\in \rad(q'_1,q'_2,q'_3 ),$  now $deq'_2=cdeF+(de)^2+ dfeg,fgq'_2=(cfg)F+(df)eg+ (fg)^2, $ and we have that $de, fg\in \rad(q'_1,q'_2,q'_3 ). $ The proof is over.

\end{example}
\begin{example} 
Consider  $S=K[a,b,c,d,e,f]$ the   ring of polynomials,  $\Bcal=\pmatrix{a&c\cr b&d\cr}$ and $F=ad-bc$, then the ideal $\Jcal=\bigcap_{i=1}^{3}\Jcal_i$, where
$$\Jcal_1=(F,e,f); \Jcal_2=(b,d,f); \Jcal_3=(a,c,e); $$
has a 2-linear resolution, 
 $\cd (\Jcal)=\projdim(S/\Jcal)=3$.
 It is not difficult to see that $3\leq \ara (\Jcal)\leq 5$, we guess that $3= \ara (\Jcal).$
\end{example}
\section{Fiber cone of codimension two lattices ideals.}
Let $\Lcal  \subset \bz^r$ be a  lattice which contains no nonnegative vectors.
Any vector $v \in \Lcal $ can be written as $v = v_+ - v_-$, where both vectors $v_+$, $v_-$ have non-negative
 coordinates. The lattice ideal $I_{\Lcal } \subset {\Scal}:= {K}[x_1,\dots ,x_r]$ is the ideal generated by all binomials
$f_v = \underline x^{v_+}- \underline x^{v_-}$, where $v$ runs in $\Lcal$. 
Prime lattice ideals are called Toric ideals, and a variety defined by a lattice ideal have a Torus action. 
By definition  the Fiber cone of $I_\mathcal L$ is the ring
$\displaystyle \Fcal(I_\Lcal) = \bigoplus_{n \ge 0} I_{\Lcal}^n/{ \bf m } I_{\Lcal}^{n+1}$: 
We have the following result from \cite{hm}
\begin{theorem}If $I _\Lcal$ is a radical ideal with $\mu $ generators and $\ K$ is infinite, then $\Fcal(I_{\Lcal})$ has dimension three, is reduced, arithmetically Cohen-Macaulay, of minimal degree. Moreover we have a description $$\Fcal(I) =  K[T_1,...,T_\mu ] \slash \widetilde{\Acal}, {\rm \ and\ }
 \widetilde{\Acal}= (\Mcal_1,\Qcal_1)\cap...\cap (\Mcal_l,\Qcal_l) , $$  where for all $k=1...,l$  $\Qcal_i$ is generated by a subset of $ \{T_1,...,T_\mu\}$ and $\Mcal_i$ is either $0$ or the ideal generated by  all the $2\times 2$ minors of a scroll matrix with only one block, which entries are subsets of $ \{T_1,...,T_\mu\}$.
\end{theorem}
\begin{corollary} If $I _\Lcal$ is a radical ideal with $\mu $ generators and $\ K$ is infinite, then $\Fcal(I_{\Lcal})$ is a set theoretical complete intersection.
\end{corollary}
The corollary is an immediate consequence of the above theorem and the theorem \ref{ara}.


\begin{thebibliography}{R-V}
\bibitem[B1]{b1} Barile, M.  {\sl Certain minimal varieties are set-theoretic complete
intersections}. Comm. Algebra {\bf 28}, (2007)1223--1239.
\bibitem[BM1]{bm1} Barile M., Morales M., {\sl On certain algebras of reduction number one,} J. Algebra {\bf 206} 
(1998), 113 -- 128.
\bibitem[BM2]{bm2} Barile M., Morales M., {\sl On the equations
    defining minimal varieties,} Comm. Alg., {\bf 28}
  (2000), 1223 -- 1239.
\bibitem[BM3]{bm3}Barile M., Morales M., {\sl On Stanley-Reisner Rings of Reduction Number One,} Ann. Sc.Nor. Sup.
 Pisa, Serie IV. Vol. {\bf XXIX} Fasc. 3. (2000), 605 -- 610.
\bibitem[BM4]{bm4} Barile M., Morales M., {\sl On unions of scrolls
    along linear spaces,} Rend. Sem. Mat. Univ. Padova, {\bf 111}
  (2004), 161 -- 178. 
 \bibitem[B]{b} Bertini, E.
{\sl Introduzione alla geometria projettiva degli iperspazi con appendice sulle curve algebriche e loro singolarità.}
Pisa: E. Spoerri. (1907).
\bibitem[BMT]{bmt} Barile, Margherita;  Morales, Marcel;  Thoma, Apostolos
{\sl On simplicial toric varieties which are set-theoretic complete intersections.} 
J. Algebra {\bf 226}, No.2, 880-892 (2000). 
\bibitem[CEP]{cep} De Concini, Corrado;  Eisenbud, David;  Procesi, Claudio
{\sl Hodge algebras.} 
Astérisque {\bf 91}, 87 p. (1982).
 \bibitem[DP]{dp} Del Pezzo, P.
{\sl Sulle superficie di ordine $n$ immerse nello spazio di $n+1$ dimensioni.} 
Nap. rend. {\bf XXIV}. 212-216. (1885).
\bibitem[EE]{ee} Eisenbud David, Evans E. Graham Jr., {\sl Every algebraic set in n-space is
the intersection of n hypersurfaces}. Invent. Math. {\bf 19}, 107-112 (1973).
\bibitem[EG]{eg} Eisenbud, David;  Goto, Shiro
{\sl Linear free resolutions and minimal multiplicity.} 
J. Algebra {\bf 88}, 89-133 (1984).
 \bibitem[EGHP]{eghp} Eisenbud D., Green M., Hulek K., Popescu S., {\sl Restricting linear syzygies: algebra and geometry,}  Compos. Math.  {\bf 141}  (2005),  no. 6, 1460--1478.
\bibitem[G]{g} Gimenez, P.; Phd thesis, University of Grenoble I.
\bibitem[GMS1]{gms1}Gimenez, P.;  Morales, M.;  Simis, A.
{\sl The analytic spread of the ideal of a monomial curve in projective 3- space}.
Eyssette, Frédéric et al., Computational algebraic geometry. 
Papers from a conference, held in Nice, France, April 21-25, 1992. Boston: Birkhäuser. Prog. Math. 109, 77-90 (1993).
\bibitem[GMS2]{gms2} Gimenez Ph., Morales M., Simis A., 
{\sl The analytical spread of the ideal of codimension 2 monomial varieties,} Result. Math. Vol {\bf 35} (1999), 250 - 259.
\bibitem[H]{h} Ha Minh Lam,   {\sl Algèbre de Rees et fibre spéciale} PhD Thesis work,
  Université J-Fourier, Grenoble, France (2006).
\bibitem[HM]{hm} Ha Minh Lam, Morales M., {\sl Fiber cone of codimension 2 lattice ideals}  To appear Comm. Alg.
\bibitem[Ha-Mo]{ha-mo}  Dao Thanh Ha, Morales Marcel, {\sl Local cohomology modules with support in 2-regular monomial ideals}  preprint, (2007).
\bibitem[Ha]{ha} Joe Harris, {\sl Algebraic Geometry, A First Course,} {\bf 1992}. Springer-Verlag, New York.
\bibitem[Mo]{m}  Morales, Marcel
 {\sl Equations des variétés monomiales en codimension deux. (Equations of monomial varieties in codimension two).}
J. Algebra {\bf 175}, No.3, 1082-1095 (1995).
\bibitem[Mo1]{mo-stci}  Morales, Marcel
 {\sl Simplicial ideals,  2-linear ideals and arithmetical rank} math.AC/0702668.
\bibitem[RV1]{rv1} Robbiano, Lorenzo;  Valla, Giuseppe
{\sl On set-theoretic complete intersections in the projective space.}
Rend. Semin. Mat. Fis. Milano {\bf 53}, 333-346 (1983).
\bibitem[S-V]{s-v} Schmitt, Th.; Vogel, W., 
{\sl Note on Set-Theoretic Intersections of Subvarieties of Projective space,}  Math. Ann.  {\bf 245} (1979), 247 - 253.
\bibitem[V]{ve} Verdi, L., 
{\sl Le curve razionali normali come intersezioni complete insiemistiche,}  Boll. UMI .  {\bf 16-A} (1979), 385--390.
\bibitem[X]{x} Xambo, S.{\sl On projective varieties of minimal degree.}
Collect. Math. {\bf 32}, 149 (1981).
\end{thebibliography}
 \end{document}